\documentclass[11pt]{article}%
\usepackage{enumerate}
\usepackage{epsfig}
\usepackage{theorem}
\usepackage{amsfonts}
\usepackage{hyperref}
\usepackage{boxedminipage}
\usepackage{color}
\usepackage{amsmath}
\usepackage{amssymb}
\usepackage{graphicx}
\usepackage{multicol}%
{\theorembodyfont{\rmfamily}

}
\newtheorem{theorem}{Theorem}

\newtheorem{corollary}[theorem]{Corollary}

\newtheorem{remark}[theorem]{Remark}

\numberwithin{equation}{section}
\numberwithin{example}{section}
\numberwithin{theorem}{section}
\begin{document}

\title{Robust bounds on risk-sensitive functionals via R\'{e}nyi divergence}
\author{Rami Atar\thanks{Research supported in part by the ISF (Grant 1349/08), the
US-Israel BSF (Grant 2008466), and the Technion fund for promotion of
research}\\Department of Electrical Engineering\\Technion--Israel Institute of Technology\\Haifa 32000, Israel
\and Kenny Chowdhary\\Sandia National Laboratories \\Livermore, CA 94550, USA
\and Paul Dupuis\thanks{Research supported in part by the Department of Energy
(DE-SC0002413,DE-SC0010539), the National Science Foundation (DMS-1317199),
and the Air Force Office of Scientific Research (FA9550-12-1-0399).}\\Division of Applied Mathematics\\Brown University\\Providence, RI 02912, USA }
\maketitle

\begin{abstract}
We extend the duality between exponential integrals and relative entropy to a
variational formula for exponential integrals involving the R\'{e}nyi
divergence. This formula characterizes the dependence of
risk-sensitive functionals and related quantities determined by tail behavior
to perturbations in the underlying distributions, in terms of the R\'{e}nyi
divergence. The characterization gives rise to upper and lower bounds that are
meaningful for all values of a large deviation scaling parameter, allowing one
to quantify in explicit terms the robustness of risk-sensitive costs. As
applications we consider problems of uncertainty quantification when aspects
of the model are not fully known, as well their use in bounding tail
properties of an intractable model in terms of a tractable one.

\vspace{\baselineskip}


\vspace{\baselineskip}


\end{abstract}

\section{Introduction}

\label{sec1}

For many models encountered in engineering, the physical sciences,
mathematical finance, and elsewhere, rare events play a key role in
determining important properties of the system. Given a system model, large
deviation theory can often be used to study the impact of rare events, and in
particular can provide both qualitative and quantitative information
\cite{frewen,dupell,demzei,shwwei}. Of course large deviation theory provides
only an asymptotic approximation, and so if non-asymptotic bounds are sought
then one can appeal to other approximations such as Monte Carlo \cite{asmgly,
blagly, dupwan5}. However, it is well known that the resulting estimates (both
asymptotic and non-asymptotic) are sensitive to the underlying assumed
distribution, owing to the fact that they are determined by tail properties of
the distributions. As a consequence, understanding the impact of modeling
errors and model uncertainty becomes especially important. Modeling
uncertainty can take many forms. For example, for some parts of the system
there may be justification for the use of distributions of a particular form,
but with parameters that are not known precisely. For other parts of the system,
however, there may not be a suitable probabilistic model, and one should
instead assume only that parameters belong to some known set.

The present paper is concerned with probabilities associated with rare events
and expected values that are largely determined by rare
events. However, the issues just raised regarding model uncertainty and
modeling error are also important for ordinary (e.g., order one)
probabilities, and expected values that are not sensitive to rare events. For
such problems, one can obtain tight bounds that hold for a well-defined family
of \textquotedblleft true\textquotedblright\ process models by computing
certain functionals with respect to a given \textquotedblleft
nominal\textquotedblright\ model, and then using the duality between
exponential integrals and relative entropy. For a detailed discussion we refer
to \cite{chodup}. Following standard terminology in the economics and control
literature, we will refer to integrals of the form
$\int_{S} e^{g}d\nu$ as \textit{risk-sensitive} functionals, where $g:S\rightarrow
\mathbb{R}$ is Borel measurable,
$S$ is a Polish space, and $\nu$ a probability measure. The well-known duality
alluded to above is
\begin{equation}
\log\int_{S}e^{g}d\nu=\sup_{\theta}\Big[\int_{S}gd\theta-R(\theta\Vert
\nu)\Big], \label{01}%
\end{equation}
where the supremum extends over all probability measures on $S$, and $R$
denotes relative entropy (see (\ref{21})). Based on this identity, the results
of \cite{chodup} give tight bounds on ordinary probabilities and expected
values, i.e., quantities of the form $\int_{S}gd\theta$. The bounds are in
terms of a maximum relative entropy distance between the nominal model, $\nu$,
and a collection of models, $\theta$, which presumably include the true model,
plus a risk-sensitive cost with respect to the nominal model. Note that the
feasibility of explicit computation, which means computing or approximating
exponential integrals, is thus linked to the choice of the nominal model.
Robust properties of controls designed on a risk-sensitive criteria were first
described in \cite{dupjampet}. By considering suitable limits such criteria
can be linked to other methods for handling model uncertainty, such as
$H^{\infty}$ control \cite{whi}.

As it turns out, the duality
\eqref{01} is not useful for bounding expectations and analyzing problems with
rare events, because the natural scaling properties are such that the
probabilities and expected values of interest should themselves be expressed
as risk-sensitive functionals (this point will be made precise later on).
However, there is a generalization of relative entropy called
\textit{R\'{e}nyi relative entropy} or \textit{R\'{e}nyi divergence}
(introduced in \cite{ren}; see Section \ref{sec2}), with which risk-sensitive
functionals can be expressed in terms of other risk-sensitive functionals.
In particular, as we shall prove, the identities
\begin{equation}
\frac{1}{\beta}\log\int_{S}e^{\beta g}d\nu=\inf_{\theta}\Big[\frac{1}{\gamma
}\log\int_{S}e^{\gamma g}d\theta+\frac{1}{\gamma-\beta}R_{\frac{\gamma}%
{\gamma-\beta}}(\nu\Vert\theta)\Big] \label{02}%
\end{equation}
and
\begin{equation}
\frac{1}{\gamma}\log\int_{S}e^{\gamma g}d\nu=\sup_{\theta}\Big[\frac{1}{\beta
}\log\int_{S}e^{\beta g}d\theta-\frac{1}{\gamma-\beta}R_{\frac{\gamma}%
{\gamma-\beta}}(\theta\Vert\nu)\Big] \label{03}%
\end{equation}
hold for any $\beta,\gamma\in\mathbb{R}\setminus\{0\}$, $\beta<\gamma$, where
for $\alpha\in\mathbb{R}\setminus\{0,1\}$ $R_{\alpha}$ denotes R\'{e}nyi
divergence of order $\alpha$ (see (\ref{10}) and \eqref{32}). Moreover,
\eqref{01} is a limit case of \eqref{03} as $\beta\rightarrow0$, with
$\gamma=1$. These identities make it possible to bound risk-sensitive
functionals with respect to the true model, $\theta$,
in terms of a risk-sensitive functional with respect to the nominal model
$\nu$. In this paper we also give elementary examples of how these bounds can be used.

As mentioned previously, one must evaluate a risk sensitive functional with
respect to a nominal model in order to turn the theoretical results into
numerical bounds. This has implications and uses that go beyond assessing
model uncertainty. In fact, it suggests an approach for bounding and
approximating rare event probabilities when evaluation of this risk-sensitive
functional is not possible or convenient for the known true model, by
replacing it with the \textquotedblleft closest\textquotedblright\ (in the
sense of R\'{e}nyi divergence) model for which the computation can be carried
out, and then bounding the R\'{e}nyi divergence between the nominal and true
models. Examples illustrating this use will be given. One can generalize to
problems of minimizing risk-sensitive costs with respect to a controlled
process, and ask for robust bounds (i.e., bounds valid for a family of process
models) in terms of the value function and optimal control for the nominal
model. This would be analogous to the robust control of order one costs by
using controls designed on the basis of risk-sensitive performance criteria
\cite{dupjampet}, and will be considered elsewhere.

We are aware of two other variational formulas for which the convex duality
relation \eqref{01} is a special case. The first is a duality formula for
$\phi$-entropy ((2.60) in \cite{mas2}, (20) in \cite{cha2}), which has
played a central role in the study of concentration inequalities
\cite{mas2}. The other is a variational formula for the $f$-divergence (a
notion similar to $\phi$-entropy), that has been used to develop
$f$-divergence estimators based on independent and identically distributed (iid) samples from each of two given
distributions. Such estimators are significant in learning problems such as
classification, dimensionality reduction, and homogeneity testing (see
\cite{nguwaijor}, \cite{rud2} for the variational formula and its uses).
Although R\'{e}nyi divergence is closely related to $f$-divergence (in
particular, the former is a certain nonlinear transformation of the latter;
see \cite{lievaj, vaj}) it seems that the representation formulas \eqref{02}
and \eqref{03} cannot be recovered from these variational characterizations.
The issue of robustness for rare events and risk-sensitive functionals has not received 
a great deal of attention. A paper that does consider the topic is \cite{kus15},
which considers the impact of varying the underlying distributions on the form of the 
large deviation rate function and related minimizers.

The rest of the paper is organized as follows. In Section \ref{sec2} we
recall the definition and some properties of R\'{e}nyi divergence, state the
variational representations based on R\'{e}nyi divergence and state some
immediate consequences. Section \ref{sec3} contains elementary applications to
functionals of empirical measures of iid outcomes, queueing, and Brownian
motion with drift, and Section \ref{sec4} concludes with the proofs of the
representation formulas.

\section{Exponential integrals and R\'{e}nyi divergence}

\label{sec2}

\subsection{Definition and properties of R\'{e}nyi divergence}

Let $({S},\mathcal{F})$ be a measurable space and let $\mathcal{P}%
=\mathcal{P}({S},\mathcal{F})$ denote the set of all probability measures on
$({S},\mathcal{F})$. We say that a measure $\mu$ on $%
({S},\mathcal{F})$ dominates $\nu\in\mathcal{P}$ if $\nu$ is absolutely
continuous with respect to $\mu$, and denote this by $\nu\ll\mu$. For two
probability measures $\nu,\theta\in\mathcal{P}$, let $\nu^\prime=\frac{d\nu
}{d\mu}$ and $\theta^\prime=\frac{d\theta}{d\mu}$ denote the Radon-Nikodym
derivatives with respect to a dominating $\sigma$-finite measure $\mu$. For
$\alpha>0$, $\alpha\neq1$, the R\'{e}nyi divergence of degree $\alpha$ of
$\nu$ from $\theta$ is defined by (cf.\ \cite{lievaj})
\begin{equation}
R_{\alpha}(\nu\Vert\theta)\doteq%
\begin{cases}
\infty & \text{if }\alpha>1\text{ and }\,\nu\not \ll \theta,\\
\displaystyle\frac{1}{\alpha(\alpha-1)}\log\int_{\{\nu^{\prime}\theta^{\prime
}>0\}}\Big(\frac{\nu^{\prime}}{\theta^{\prime}}\Big)^{\alpha}d\theta &
\text{otherwise.}%
\end{cases}
\label{10}%
\end{equation}
We follow \cite{lievaj} in defining $R_\alpha$ with the factor $\frac
{1}{\alpha(\alpha-1)}$ rather than $\frac{1}{\alpha-1}$, which is also a
common choice \cite{bha,ren,vaj}. When $\nu$ and $\theta$ are mutually
absolutely continuous, this expression can be written without reference to a
dominating measure, namely
\[
R_{\alpha}(\nu\Vert\theta)=\frac{1}{\alpha(\alpha-1)}\log\int_{S}%
\Big(\frac{d\nu}{d\theta}\Big)^{\alpha}d\theta=\frac{1}{\alpha(\alpha-1)}%
\log\int_{S}\Big(\frac{d\theta}{d\nu}\Big)^{1-\alpha}d\nu.
\]
The definition of $R_\alpha$ is extended to $\alpha=1$ by letting $R_1=R$ be
the relative entropy, or the Kullback-Liebler divergence, defined by
\begin{equation}
R(\nu\Vert\theta)\doteq%
\begin{cases}
\infty & \text{if }\nu\not \ll \theta,\\
\displaystyle\int_{\{\nu^{\prime}\theta^{\prime}>0\}}\frac{\nu^{\prime}%
}{\theta^{\prime}}\log\frac{\nu^{\prime}}{\theta^{\prime}}d\theta &
\text{otherwise.}%
\end{cases}
\label{21}%
\end{equation}
The definitions do not depend on the choice of the dominating measure, and
since $\nu+\theta$ automatically dominates $\nu$ and $\theta$, $R_\alpha
(\nu\Vert\theta)$ is well defined for all pairs $(\nu,\theta)\in
\mathcal{P}^2$. For a proof of independence from the dominating measure as
well as various properties of $R_\alpha$, see
\cite{golpasyar,lievaj,vaj,ervhar}. To mention a few of these properties, let
$\nu$ and $\theta$ be fixed. Then $\alpha\mapsto\alpha R_\alpha(\nu\Vert
\theta)$ is nondecreasing as a map from $(0,\infty)$ to $[0,\infty]$, and
continuous from the left (thus $R=\lim_{\alpha\uparrow 1}R_\alpha$).
If $\nu$ and $\theta$ are mutually singular then $R_\alpha(\nu\Vert\theta)$ is
infinite everywhere. Otherwise, it is finite and continuous on $(0,\bar
{\alpha})$,
where $\bar{\alpha}\doteq\sup\{\alpha:R_\alpha(\nu\Vert\theta)<\infty\}\geq1$.
Moreover, for every $\alpha>0$, $R_\alpha(\nu\Vert\theta)=0$ if and only if $\nu=\theta$.

A further useful property is the identity $R_{\alpha}(\nu\Vert\theta
)=R_{1-\alpha}(\theta\Vert\nu)$, which holds for every $\alpha\in(0,1)$. We
will use it to extend the definition of $R_{\alpha}$ to $\alpha\in
\mathbb{R}\setminus\{0,1\}$. Namely, we set
\begin{equation}
R_{\alpha}(\nu\Vert\theta)\doteq R_{1-\alpha}(\theta\Vert\nu),\qquad\alpha<0.
\label{32}%
\end{equation}
This definition is consistent with the definition of $R_{\alpha}$, $\alpha
\in\mathbb{R}$, given in (2.10) of \cite{lievaj}, as follows from Remark 2.13
of \cite{lievaj}.

\subsection{Variational representations for exponential integrals}

The variational representation for exponential integrals \eqref{01} is very
closely related to the theory of large deviations, and in fact can serve as
the natural starting point for the large deviations analysis of any system
\cite{dupell4}. It also gives an inequality that allows for robust bounds on
ordinary costs with respect to a \textquotedblleft true\textquotedblright%
\ measure in terms of risk-sensitive costs for a \textquotedblleft
nominal\textquotedblright\ model plus relative entropy distance between the
two. However, as noted in the Introduction, this variational representation
does not seem to be useful when bounding risk-sensitive costs. The variational
representations in Theorem \ref{th1} give useful bounds in that respect. A
particular case of (\ref{11-}) appears in \cite{dvitod}. The proof of the
theorem is given in Section \ref{sec4}.

\begin{theorem}
\label{th1} Let $\beta$ and $\gamma$ be members of $\mathbb{R}\setminus\{0\}$,
with $\beta<\gamma$. Let $\nu\in\mathcal{P}$. Then for any bounded and
measurable $g:S\rightarrow\mathbb{R}$, one has
\begin{equation}
\frac{1}{\beta}\log\int_{S}e^{\beta g}d\nu=\inf_{\theta\in\mathcal{P}%
}\Big[\frac{1}{\gamma}\log\int_{S}e^{\gamma g}d\theta+\frac{1}{\gamma-\beta
}R_{\frac{\gamma}{\gamma-\beta}}(\nu\Vert\theta)\Big], \label{11-}%
\end{equation}
where the infimum is uniquely attained at $d\theta=e^{-(\gamma-\beta)g}d\nu/Z
$, $Z=\int_{S}e^{-(\gamma-\beta)g}d\nu$. In addition,
\begin{equation}
\frac{1}{\gamma}\log\int_{S}e^{\gamma g}d\nu=\sup_{\theta\in\mathcal{P}%
}\Big[\frac{1}{\beta}\log\int_{S}e^{\beta g}d\theta-\frac{1}{\gamma-\beta
}R_{\frac{\gamma}{\gamma-\beta}}(\theta\Vert\nu)\Big], \label{29-}%
\end{equation}
where the supremum is uniquely attained at $d\theta=e^{(\gamma-\beta)g}d\nu/Z
$, $Z=\int_{S}e^{(\gamma-\beta)g}d\nu$.
\end{theorem}

\begin{remark}
\label{rem0} \emph{Setting }$\beta=\alpha-1$\emph{ and }$\gamma=\alpha$\emph{
gives}
\begin{equation}
\frac{1}{\alpha-1}\log\int_{S}e^{(\alpha-1)g}d\nu=\inf_{\theta\in\mathcal{P}%
}\Big[\frac{1}{\alpha}\log\int_{S}e^{\alpha g}d\theta+R_{\alpha}(\nu
\Vert\theta)\Big],\label{11}%
\end{equation}
\emph{and}
\begin{equation}
\frac{1}{\alpha}\log\int_{S}e^{\alpha g}d\nu=\sup_{\theta\in\mathcal{P}%
}\Big[\frac{1}{\alpha-1}\log\int_{S}e^{(\alpha-1)g}d\theta-R_{\alpha}%
(\theta\Vert\nu)\Big],\label{29}%
\end{equation}
\emph{for all }$\alpha\neq0$\emph{, }$\alpha\neq1$\emph{. Although \eqref{11}
and \eqref{29} are special cases of \eqref{11-} and \eqref{29-}, the latter
can be recovered from the former (in fact with }$\alpha$\emph{ in the range
}$\alpha>0$\emph{, }$\alpha\neq1$\emph{), as shown in the proof of Theorem
\ref{th1}. In most of what follows we will work with \eqref{11} and
\eqref{29}}.
\end{remark}

\begin{remark}
\emph{By taking the formal limit }$\alpha\rightarrow1$\emph{, we obtain from
(\ref{11}) the identity }%
\[
\int_{S}gd\nu=\inf_{\theta\in\mathcal{P}}\Big[\log\int_{S}e^{g}d\theta
+R(\nu\Vert\theta)\Big]
\]
\emph{and from (\ref{29}) the well-known convex duality formula (see
\cite{chodup,dupell4,dupjampet})}
\[
\log\int_{S}e^{g}d\nu=\sup_{\theta\in\mathcal{P}}\Big[\int_{S}gd\theta
-R(\theta\Vert\nu)\Big].
\]
\emph{Note that one can also take }$\alpha\rightarrow0$\emph{\ in \eqref{11}
and \eqref{29}, in which case }$\alpha R_{\alpha}(\nu\Vert\theta
)\rightarrow-\log\theta(\nu^{\prime}>0)$\emph{, recovering the simple fact}
\[
0=\inf_{\theta\in\mathcal{P}}[-\log\theta(\nu^{\prime}>0)]=\sup_{\theta
\in\mathcal{P}}\log\nu(\theta^{\prime}>0).
\]

\end{remark}

The main purpose of this paper is to observe the following inequalities that
follow from \eqref{11} and \eqref{29}, and to discuss how they can be used to
study robustness of risk-sensitive functionals.

\begin{corollary}
\label{cor:1}Assume $\alpha>1$, $\theta\in\mathcal{P}$, $\nu\in\mathcal{P}$,
and let $g:S\rightarrow\mathbb{R}$ be any measurable function. Then
\begin{align}
\frac{1}{\alpha-2}\log\int_{S}e^{(\alpha-2)g}d\nu-R_{\alpha-1}(\nu
\Vert\theta) &  \leq\frac{1}{\alpha-1}\log\int_{S}e^{(\alpha-1)g}d\theta\label{13}\\
&  \leq\frac{1}{\alpha}\log\int_{S}e^{\alpha g}d\nu+R_{\alpha}(\theta
\Vert\nu),\nonumber
\end{align}
where the first inequality also requires $\alpha>2$. Also, on the left hand
side of (\ref{13}) we interpret $\infty-\infty$ as $-\infty$.
\end{corollary}

See Section \ref{sec4} for the proof. Similar inequalities can be deduced when
$\alpha\in(0,1)$, but for our present purposes they do not seem to be
particularly useful.

The following interpretation of Corollary \ref{cor:1} will be useful in the
examples presented in the next section. By considering $g=0$ on $A\in
\mathcal{F}$ and $g=-M$ on $A^{c}$, and then sending $M\rightarrow\infty$, one
obtains that for any event $A$%
\begin{equation}
\frac{1}{\alpha-2}\log\nu(A)-R_{\alpha-1}(\nu\Vert\theta)\leq\frac
{1}{\alpha-1}\log\theta(A)\leq\frac{1}{\alpha}\log\nu(A)+R_{\alpha}(\theta
\Vert\nu), \label{31}%
\end{equation}
with the same restrictions on $\alpha$ as in the corollary.

\section{Elementary applications}

\label{sec3}In this section we show how Corollary \ref{cor:1} can be used to
provide robust bounds of the sort described in the Introduction. The examples
are intended only to illustrate the main ideas, and limited to problems where
the driving noises are distributed according to product measure. When
assessing probabilities and expected values associated with rare events, it is
important to keep in mind that it is usually \textit{relative errors}, and not
absolute errors, that are important. Also, it is generally the case that
approximations are of an asymptotic nature as some scaling parameter tends to
a limit. For light-tailed processes, the scaling is exponential in the
parameter. As we will see, this fits in very nicely with the form of the
inequalities in (\ref{13}). 

As described in the introduction, one should have
in mind two scenarios. In one case, we think of $\theta$ as a probability measure
of interest for which the large deviation functional may be hard to compute,
and of $\nu$ as an alternative that is more tractable. In the other case we
are not sure of the model, with the nominal model $\nu$ a sort of
\textquotedblleft best guess\textquotedblright\ and $\theta$ the true model.

\subsection{Functionals of the empirical measure}

Suppose that $S={\mathbb{R}}^{n}$, where $n$ is the scaling parameter. Let
$\theta^{n}$ and $\nu^{n}$ be product probability measures on $S$, with
marginals $\theta_{i}^{n}$ and $\nu_{i}^{n}$. Assume $\nu_{i}^{n}%
=\nu_{1}$, so the nominal model corresponds to an iid sequence. Then
(cf.\ \cite{golpasyar})
\[
{\Delta}_{\alpha}^{n}\doteq\frac{1}{n}R_{\alpha}(\theta^{n}\Vert\nu^{n}%
)=\frac{1}{n}\sum_{i=1}^{n}R_{\alpha}(\theta_{i}^{n}\Vert\nu_{i}^{n})=\frac
{1}{n}\sum_{i=1}^{n}R_{\alpha}(\theta_{i}^{n}\Vert\nu_{1}).
\]
Let $X_{n}$ denote the canonical process. If the $X_{n}$ are also iid under
$\theta^{n}$ with marginal $\theta_{1}$, then ${\Delta}_{\alpha}^{n}=R_{\alpha}%
(\theta_{1}\Vert\nu_{1})$ for every $n$. Consider the empirical measure
$L_{n}=n^{-1}\sum_{i=1}^{n}\delta_{X_{i}}$ as a random element of the space
$\mathcal{P}_{{\mathbb{R}}}=\mathcal{P}({\mathbb{R}},\mathcal{R})$, equipped
with the topology of weak convergence, and fix any measurable function
$G:\mathcal{P}\rightarrow{\mathbb{R}}$. Then with ${\mathbb{E}}_{\theta}$ and
${\mathbb{E}}_{\nu}$ denoting expectation with respect to the indicated
distribution, we can take $g(X_{n})=nG(L_{n})$ in Corollary \ref{cor:1} to
get
\begin{equation}
\frac{1}{n}\frac{1}{\alpha-1}\log{\mathbb{E}}_{\theta}e^{(\alpha-1)nG(L_{n})}%
\leq\frac{1}{n}\frac{1}{\alpha}\log{\mathbb{E}}_{\nu}e^{\alpha nG(L_{n}%
)}+R_{\alpha}(\theta_{1}\Vert\nu_{1}).\label{16}%
\end{equation}
(and also if desired a corresponding lower bound).

If $G$ is continuous and $\theta$ corresponds to an iid sequence, then in this
very simple setting one could use Sanov's theorem to evaluate the limit
behavior of the two terms, and obtain
\begin{equation}
\frac{1}{\alpha-1}\sup_{\lambda\in\mathcal{P}_{{\mathbb{R}}}}[(\alpha
-1)G(\lambda)-R(\lambda\Vert\theta_{1})]\leq\frac{1}{\alpha}\sup_{\lambda
\in\mathcal{P}_{{\mathbb{R}}}}[\alpha G(\lambda)-R(\lambda\Vert\nu
_{1})]+R_{\alpha}(\theta_{1}\Vert\nu_{1}). \label{14}%
\end{equation}
The strength of the general inequalities based on R\'{e}nyi divergence is that
the bound \eqref{16} holds for all $n$, and moreover does not require that
$\theta$ correspond to an iid sequence.

We can make \eqref{16} and \eqref{14} more concrete by considering, for example, Gaussian
distributions $\theta_{1}=\mathcal{N}(\mu_{1},\sigma_{1}^{2})$ and $\nu
_{1}=\mathcal{N}(\mu_{2},\sigma_{2}^{2})$. In this case
\begin{equation}
R_{\alpha}(\theta_{1}\Vert\nu_{1})=%
\begin{cases}
\frac{1}{\alpha}\log\frac{\sigma_{2}}{\sigma_{1}}+\frac{1}{2\alpha(\alpha
-1)}\log\frac{\sigma_{2}^{2}}{\sigma_{\alpha}^{2}}+\frac{1}{2\alpha}%
\frac{\alpha(\mu_{1}-\mu_{2})^{2}}{\sigma_{\alpha}^{2}}, & \text{ if }%
\sigma_{\alpha}^{2}\doteq\alpha\sigma_{2}^{2}+(1-\alpha)\sigma_{1}^{2}>0,\\
\infty, & \text{ otherwise.}%
\end{cases}
\label{17}%
\end{equation}
If $G(L_{n})=c\langle1,L_{n}\rangle=n^{-1}c(X_{1}+\cdots+X_{n})$ for some
constant $c$ and $\nu_{1}=\mathcal{N}(0,1)$, then ${\mathbb{E}}_{\nu
}e^{\alpha nG(L_{n})}={\mathbb{E}}_{\nu}e^{\alpha c(X_{1}+\cdots+X_{n}%
)}=e^{\alpha^{2}c^{2}n/2}$ and \eqref{16} says that for every $\theta$ under
which $X_{n}$ are iid,
\begin{equation}
\frac{1}{n}\frac{1}{\alpha-1}\log{\mathbb{E}}_{\theta}e^{(\alpha-1)c(X_{1}%
+\cdots+X_{n})}\leq R_{\alpha}(\theta_{1}\Vert\nu_{1})+\frac{\alpha c^{2}}%
{2}.\label{18}%
\end{equation}
In (\ref{18}) one obtains equality if $\theta$ is $\mathcal{N}(c,1)$, as can be
verified using \eqref{17}. As a result, \eqref{18} is tight in the following
sense. Fix $\alpha>1$ and a constant $d>0$. Consider the family of $\theta_{1}$
for which $R_{\alpha}(\theta_{1}\Vert\nu_{1})\leq d$. With this notation,
\eqref{18} states that for $\theta_{1}$ in this family,
\[
\frac{1}{n}\frac{1}{\alpha-1}\log{\mathbb{E}}_{\theta}e^{(\alpha-1)c(X_{1}%
+\cdots+X_{n})}\leq d+\frac{\alpha c^{2}}{2}.
\]
Moreover one can find $c$ and a $\theta_{1}$ in the family such that this display
holds with equality. Indeed, $c$ is chosen so that $\frac{1}{2}(\alpha
-1)c^{2}=d$ (namely, $c=\pm\sqrt{2d/(\alpha-1)}$) and $\theta_{1}=\mathcal{N}%
(c,1)$.

\subsection{A sample path large deviation example}

We next discuss a well-known example from queueing analysis. Lindley's
recursion
\[%
\begin{cases}
Q_{n}=(Q_{n-1}+X_{n}-C)^{+}, & n\geq1,\\
Q_{0}=0, &
\end{cases}
\]
describes the queue length $Q_{n}$ in an initially empty queueing system where
$X_{n}\geq0$ arrivals occur at time $n\geq1$, and the server is capable of
serving $C$ customers at each time slot. Denoting $S_{0}=0$ and $S_{n}%
=X_{1}+\cdots+X_{n}$, the solution to this recursion is given by
\[
Q_{n}=\max_{0\leq i\leq n}[S_{n}-S_{i}-Ci].
\]
Assume that the system is stable in the sense that $C>{\mathbb{E}}_{\nu
}[X_{1}]=1$. Consider the space-time rescaled processes $\bar{S}^{n}%
(t)=n^{-1}S_{[nt]}$ and $\bar{Q}^{n}(t)=n^{-1}Q_{[nt]}$, $t\geq0$, and given a
constant $b>0$, let the \textit{buffer overflow} event be given by
\[
A_{n}=\left\{  \max_{t\in\lbrack0,1]}\bar{Q}^{n}(t)>b\right\}  .
\]

The large deviation asymptotic behavior of this sequence of events has been
studied in general; see, for example, \cite{ana2}, and Section 11.7 of
\cite{shwwei}. Here we will focus on a simple special case. Assume that under
$\nu$, $X_{n}$ are iid standard Poisson. Let $\mathcal{AC}([0,1]:\mathbb{R})$ [resp., $\mathcal{D}([0,1]:\mathbb{R})$] denote the
space of functions that are absolutely continuous
[resp., right continuous with left limits] and that map $[0,1]$ to $\mathbb{R}$.
Equip $\mathcal{D}([0,1]:\mathbb{R})$ with the Skorohod $J_1$ topology.
The processes $\bar{S}^{n}$ are known to satisfy a sample-path large deviation
principle in $\mathcal{D}([0,1]:\mathbb{R})$
with the rate
function $I$ given by
\[
I(\varphi)=%
\begin{cases}
\int_{0}^{1}\ell(\dot{\varphi}(t))dt, & \text{if }\varphi\in\mathcal{AC}([0,1]:\mathbb{R})%
,\varphi(0)=0,\\
\infty, & \text{otherwise,}%
\end{cases}
\]
where, with the convention $0\log0=0$,
\[
\ell(x)=%
\begin{cases}
x\log x-x+1, & \text{if }x\geq0,\\
\infty, & \text{if }x<0,
\end{cases}
\]
\cite[Theorem 6.1(b)]{puhwhi}. Hence $\lim_{n\rightarrow\infty}\frac{1}%
{n}\log\nu(A_n)=-c$, where
\[
c\doteq\inf\Big\{\int_{0}^{1}\ell(\dot{\varphi}(t))dt:\varphi\in
\mathcal{AC},\varphi(0)=0,\max_{0\leq s\leq t\leq1}\varphi(t)-\varphi
(s)-C(t-s)\geq b\Big\}
\]
can be found explicitly. Let $m^\ast$ and $t^\ast$ denote the minimum of
$t\ell(C+\frac{b}{t})$ over $t>0$ and the unique minimizer, respectively.
Then
\[
c=%
\begin{cases}
\ell(C+b) & \text{if }t^{\ast}\geq1,\\
m^{\ast} & \text{if }t^{\ast}<1.
\end{cases}
\]
Note that the event $A_n$ depends only on $X_1,\ldots,X_n$. If $\theta$ is any
probability measure under which $X_n$ are iid and $R_\alpha(\theta_1\Vert
\nu_1)\leq d_1$ and $R_{\alpha-1}(\nu_1\Vert\theta_1)\leq d_2$ for constants
$d_1,d_2$, then we obtain from \eqref{31} that for all $n$%
\[
\frac{1}{n}\frac{1}{\alpha-2}\log{\mathbb{P}}_{\nu}(A_{n})-d_{2}\leq
\frac{1}{n}\frac{1}{\alpha-1}\log{\mathbb{P}}_{\theta}(A_{n})\leq\frac{1}{n}%
\frac{1}{\alpha}\log{\mathbb{P}}_{\nu}(A_{n})+d_{1},
\]
or%
\[
{\mathbb{P}}_{\nu}(A_{n})^{\frac{\alpha-1}{\alpha-2}}e^{n(\alpha-1)d_{2}%
}\leq{\mathbb{P}}_{\theta}(A_{n})\leq{\mathbb{P}}_{\nu}(A_{n})^{\frac
{\alpha-1}{\alpha}}e^{n(\alpha-1)d_{1}}.
\]
In particular,
\[
-\frac{\alpha-1}{\alpha-2}c-(\alpha-1)d_{2}+o(1)\leq\frac{1}{n}\log
\mathbb{P}_{\theta}(A_{n})\leq-\frac{\alpha-1}{\alpha}c+(\alpha-1)d_{1}+o(1),
\]
as $n\rightarrow\infty$. More generally, the same conclusions hold if $\theta$ is
any product measure under which
\[
\frac{1}{n}\sum_{i=1}^{n}R_{\alpha}(\theta_{i}^{n}\Vert\nu_{1})\leq
d_{1},\quad\frac{1}{n}\sum_{i=1}^{n}R_{\alpha-1}(\nu_{1}\Vert\theta_{i}%
^{n})\leq d_{2}.
\]

\subsection{Brownian motion with drift}

Let $B_{t}$ be standard Brownian motion on $0\leq t\leq1$ and let $P$ be the
corresponding standard Wiener measure on $\mathcal{C}([0,1]:\mathbb{R})$. Let $Q$ be the
measure induced by Brownian motion with constant drift, i.e.,
\[
X_{t}=B_{t}+\mu t,
\]
where $\mu\in\mathbb{R}$. Also, let $\tilde{Q}$ be the measure induced by the
paths of the solution $X$ to the stochastic differential equation (SDE)
\[
dX_{t}=m(X_{t})dt+dB_{t},\qquad X_{0}=0,
\]
for measurable $m$, where, by assumption, weak existence and uniqueness hold.
A simple calculation based on Girsanov's theorem yields that the R\'{e}nyi
divergence between $Q$ and $P$ is given by
\begin{equation}
R_{\alpha}(Q\Vert P)=R_{\alpha}(P\Vert Q)=\frac{\mu^{2}}{2}, \label{eq:qp}%
\end{equation}
and that, if $|m(x)|\leq|\mu|$ for all $x$, then
\begin{equation}
R_{\alpha}(\tilde{Q}\Vert P)\leq\frac{\mu^{2}}{2},\qquad R_{\alpha}%
(P\Vert\tilde{Q})\leq\frac{\mu^{2}}{2}. \label{04}%
\end{equation}

Let $A$ be the event that the path exceeds a certain level $K>0$:
\[
A\doteq\left\{  \omega:\sup_{0\leq t\leq1}X_{t}>K\right\}  .
\]
The exceedance probability under the measure $Q$, which represents the
probability of Brownian motion with constant drift exceeding $K$, is given
(see \cite[\S 2.1]{borsal}) by
\[
Q(A)=\frac{1}{2}\text{erfc}\left(  \frac{K-\mu}{\sqrt{2}}\right)  +\frac{1}%
{2}e^{2\mu K}\text{erfc}\left(  \frac{K+\mu}{\sqrt{2}}\right)  ,
\]
where $\text{erfc}(x)=\frac{2}{\sqrt{\pi}}\int_{x}^{\infty}e^{-v^{2}}dv$, and
under standard Wiener measure,
\begin{equation}
P(A)=\sqrt{\frac{2}{\pi}}\int_{K}^{\infty}e^{-x^{2}/2}dx=\text{erfc}\left(
\frac{K}{\sqrt{2}}\right)  . \label{eq:evalp}%
\end{equation}
We would like to identify the bounds on $Q(A)$ and $\tilde{Q}(A)$ that
Corollary \ref{cor:1} provides. In particular, by \eqref{31}
\[
\frac{1}{\alpha-2}\log P(A)-R_{\alpha-1}(P\Vert Q)\leq\frac{1}{\alpha-1}\log
Q(A)\leq\frac{1}{\alpha}\log P(A)+R_{\alpha}(Q\Vert P),
\]
where the right hand side is valid for $\alpha>1$ and the left hand side is
valid for $\alpha>2$. By (\ref{eq:qp}) and (\ref{eq:evalp}) this gives
\[
\frac{1}{\alpha-2}\log\text{erfc}\left(  \frac{K}{\sqrt{2}}\right)  -\frac
{\mu^{2}}{2}\leq\frac{1}{\alpha-1}\log Q(A)\leq\frac{1}{\alpha}\log
\text{erfc}\left(  \frac{-K}{\sqrt{2}}\right)  +\frac{\mu^{2}}{2},
\]
or in probability scale
\[
\text{erfc}\left(  \frac{-K}{\sqrt{2}}\right)  ^{\frac{\alpha-1}{\alpha-2}%
}e^{-(\alpha-1)\frac{\mu^{2}}{2}}\leq Q(A)\leq\text{erfc}\left(  \frac
{-K}{\sqrt{2}}\right)  ^{\frac{\alpha-1}{\alpha}}e^{(\alpha-1)\frac{\mu^{2}%
}{2}}.
\]
By \eqref{04}, the same conclusion holds for $\tilde{Q}(A)$. 

\begin{figure}
[ptb]
\begin{center}
\includegraphics[
height=3.0242in,
width=4.0292in
]%
{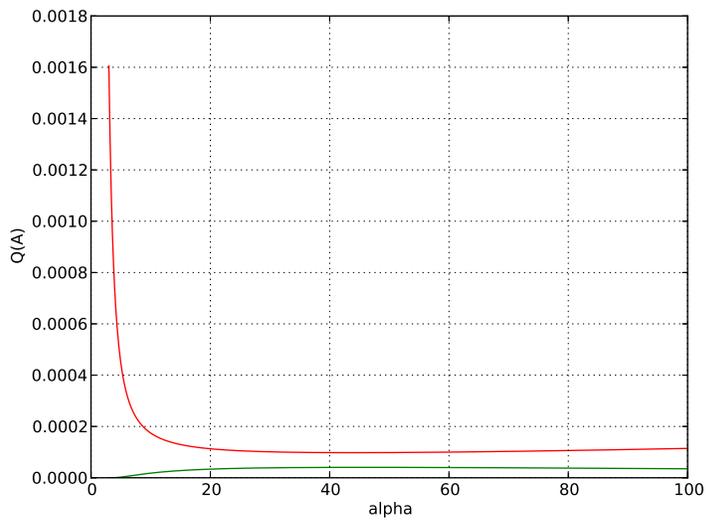}%
\caption{Upper and lower bounds for $Q(A) $ in probability scale.}%
\label{prob_scale2}%
\end{center}
\end{figure}

\begin{figure}
[ptb]
\begin{center}
\includegraphics[
height=3.0242in,
width=4.0292in
]%
{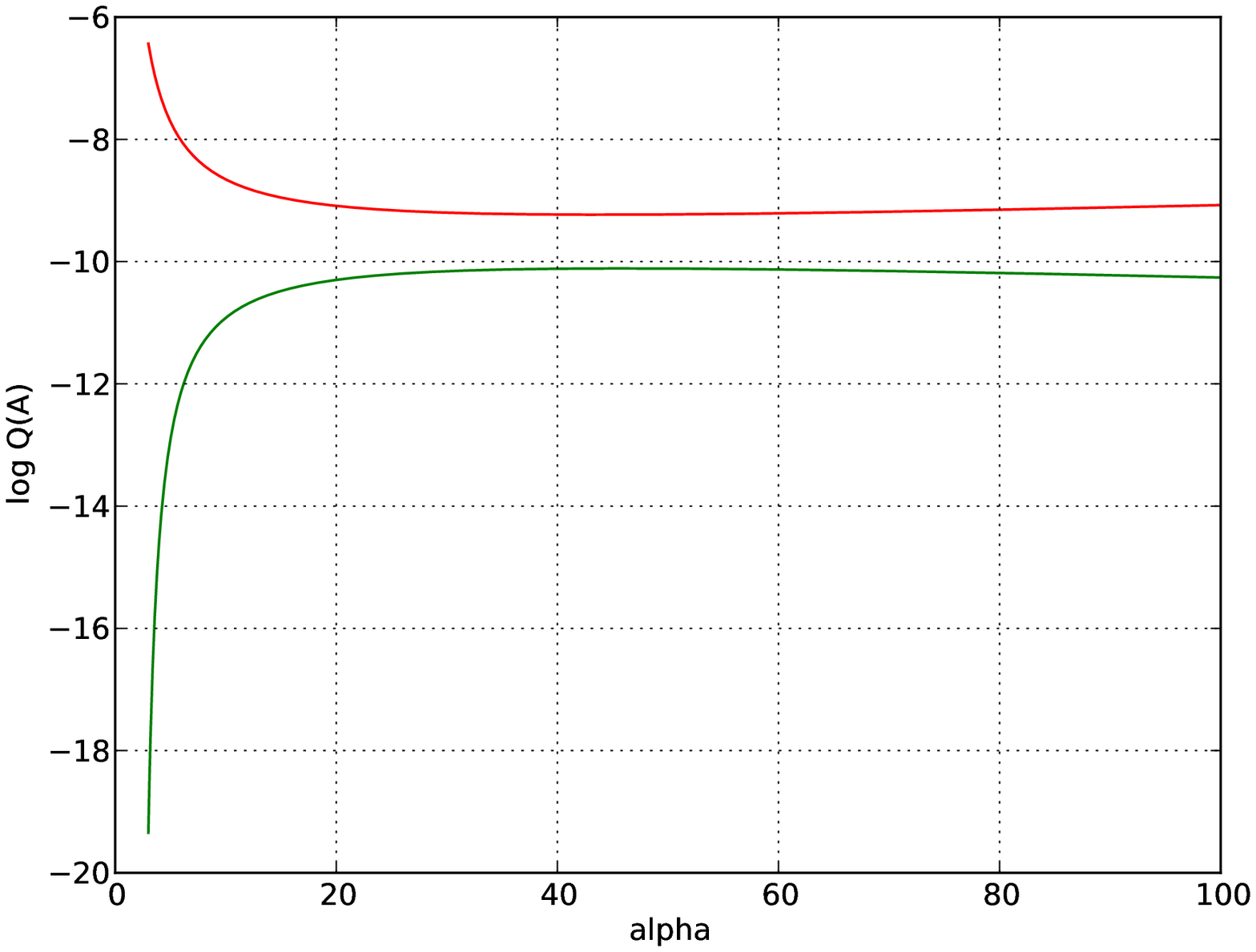}%
\caption{Upper and lower bounds for $Q(A)$ in log probability scale.}%
\label{logprob_scale2}%
\end{center}
\end{figure}

To illustrate these upper and lower bounds, we consider Brownian motion with
constant drift with $|\mu|\leq.1$ so that $R_{\alpha}(P\Vert Q)\leq.005$ and
$K=4$. Note that with $K=4$,
\[
P(A)\approx6.33\times10^{-5}.
\]
Figures \ref{prob_scale2} and \ref{logprob_scale2} show the upper bounds in
probability and log-probability scale, respectively, plotted as a function of
$\alpha\in\lbrack3,100]$.%


As another example involving the measures $P$, $Q$ and $\tilde{Q}$, consider
 the random variable
\[
H(t)=\inf\left\{  s\in\lbrack0,t]:X_{s}=\sup_{u\in\lbrack0,t]}X_{u}\right\}
,\qquad t\geq0.
\]
The Laplace transform of $H(t)$ in the case of the standard Wiener measure is given by
\[
\mathbb{E}_{P}[e^{-\gamma H(t)}]=e^{-\gamma t/2}I_{0}\left(  \frac{\gamma
t}{2}\right)  .
\]
For the case of constant drift,
\[
\mathbb{E}_{Q}[e^{-\gamma H(t)}]=\left(  \frac{e^{-\gamma t-\mu^{2}t/2}}%
{\sqrt{\pi t}}+\frac{\mu e^{-\gamma t}}{\sqrt{2}}\text{erfc}\left(  -\frac
{\mu\sqrt{t}}{\sqrt{2}}\right)  \right)  \ast\left(  \frac{e^{-\mu^{2}t/2}%
}{\sqrt{\pi t}}-\frac{\mu}{\sqrt{2}}\text{erfc}\left(  \frac{\mu\sqrt{t}%
}{\sqrt{2}}\right)  \right)
\]
where $f(t)\ast g(t)$ denotes the convolution of $f$ and $g$ evaluated at $t$
(see \cite{borsal}). There is no explicit expression for the case of a SDE.
To obtain bounds on the behavior under $Q$ and $\tilde{Q}$ we apply Corollary
\ref{cor:1}, which gives
\begin{align*}
&  \frac{1}{\alpha-2}\left[  -\frac{(\alpha-2)\gamma t}{2}+\log I_{0}\left(
\frac{(\alpha-2)\gamma t}{2}\right)  \right]  -\frac{\mu^{2}t}{2}\\
&  \hspace{8em}\leq\frac{1}{\alpha-1}\log\mathbb{E}_{Q}[e^{(\alpha-1)\gamma
H(t)}]\\
&  \hspace{8em}\leq\frac{1}{\alpha}\left[  -\frac{\alpha\gamma t}{2}+\log
I_{0}\left(  \frac{\alpha\gamma t}{2}\right)  \right]  +\frac{\mu^{2}t}{2}.
\end{align*}
As before, the same upper and lower bounds are valid for $\tilde{Q}$ as well.

\section{Proofs of \textbf{Theorem \ref{th1} and Corollary \ref{cor:1}}}

\label{sec4}

\textbf{Proof of Theorem \ref{th1}.} The main part of the proof will be to
show the validity of \eqref{11} and \eqref{29} for all $\alpha>0$, $\alpha
\ne1$. Before proving these identities, let us show that they imply
\eqref{11-} and \eqref{29-}. First, note that \eqref{11} and \eqref{29} for
$\alpha>0$, $\alpha\ne1$ imply \eqref{11} and \eqref{29} for all $\alpha
\in\mathbb{R}\setminus\{0,1\}$. Indeed, if $\alpha<0$ then \eqref{11} with
$\bar\alpha=1-\alpha>1$ and $\bar g= -g$ reads
\[
\frac{1}{\bar\alpha-1}\log\int e^{(\bar\alpha-1)\bar g}d\theta=\inf_{\theta
\in\mathcal{P}}\Big[\frac{1}{\bar\alpha}\log\int e^{\bar\alpha\bar g}%
d\theta+R_{\bar\alpha}(\nu\|\theta)\Big].
\]
Expressed in terms of $\alpha$ and $g$,
\[
-\frac{1}{\alpha}\log\int e^{\alpha g}d\theta=\inf_{\theta\in\mathcal{P}%
}\Big[-\frac{1}{\alpha-1}\log\int e^{(\alpha-1)g}d\theta+R_{\alpha}%
(\theta\|\nu)\Big],
\]
where we used \eqref{32}. Multiplying by $(-1)$ establishes the validity of
\eqref{29} for $\alpha<0$. In a similar way, the validity of \eqref{11} for
$\alpha<0$ follows from that of \eqref{29} for $\bar\alpha>1$.

Next, to show that \eqref{11} and \eqref{29} with $\alpha\in\mathbb{R}%
\setminus\{0,1\}$ imply \eqref{11-} and \eqref{29-}, fix $\beta$ and $\gamma$
in $\mathbb{R}\setminus\{0\}$, $\beta<\gamma$. Apply \eqref{11} with
$\alpha=\frac{\gamma}{\gamma-\beta}$ and $g=(\gamma-\beta)f$ (note that
$\alpha\notin\{0,1\}$) and divide by $\gamma-\beta$ to get \eqref{11-} (with
$f$ in place of $g$). In a similar way, \eqref{29-} follows from \eqref{29}.

We turn to proving (\ref{11}) for $\alpha>0$, $\alpha\neq1$. Fix $\nu$, and
consider first the case $\alpha>1$. Given any $\theta$, let $\mu=\mu(\theta)$
be a measure dominating both $\nu$ and $\theta$, and denote by $\nu^{\prime}$
and $\theta^{\prime}$ the corresponding densities. Define $\lambda
\in\mathcal{P}$ by $d\lambda=e^{-g}d\nu/Z$ where $Z=\int_{S}e^{-g}d\nu$, and
let $\lambda^{\prime}$ be the density of $\lambda$ with respect to $\mu$.
First suppose that $\theta$ dominates $\nu$. Then $\lambda^{\prime}/\nu^{\prime
}=e^{-g}/Z$, and so
\begin{align}
\log\int_{S}e^{\alpha g}d\theta &  \geq\log\int_{\{\nu^{\prime}>0\}}e^{\alpha
g}d\theta\label{24}\\
&  =\log\int_{\{\nu^{\prime}>0\}}\frac{1}{Z}\frac{\nu^{\prime}}{\lambda
^{\prime}}e^{(\alpha-1)g}d\theta\nonumber\\
&  =\log\int_{S}\frac{1}{Z}\frac{\nu^{\prime}\theta^{\prime}}{\lambda^{\prime
}}e^{(\alpha-1)g}d\mu\nonumber\\
&  =\log\int_{S}\frac{1}{Z}\frac{\theta^{\prime}}{\lambda^{\prime}}%
e^{(\alpha-1)g}d\nu.\nonumber
\end{align}
Moreover, since $\nu\ll\theta$ $\mu\{\nu^{\prime}\theta^{\prime}>0\}=\mu
\{\nu^{\prime}>0\}$, and therefore
\begin{align*}
R_{\alpha}(\nu\Vert\theta) &  =\frac{1}{\alpha(\alpha-1)}\log\int
_{\{\nu^{\prime}\theta^{\prime}>0\}}\Big(\frac{\nu^{\prime}}{\theta^{\prime}%
}\Big)^{\alpha}d\theta\\
&  =\frac{1}{\alpha(\alpha-1)}\log\int_{\{\nu^{\prime}>0\}}\Big(\frac
{\nu^{\prime}}{\theta^{\prime}}\Big)^{\alpha-1}d\nu\\
&  =\frac{1}{\alpha(\alpha-1)}\log\int_{\{\nu^{\prime}>0\}}Z^{\alpha
-1}\Big(\frac{\lambda^{\prime}}{\theta^{\prime}}\Big)^{\alpha-1}%
e^{(\alpha-1)g}d\nu.
\end{align*}
Thus with $d\tilde{\nu}=e^{(\alpha-1)g}d\nu$,
\begin{equation}
\frac{1}{\alpha}\log\int_{S}e^{\alpha g}d\theta+R_{\alpha}(\nu\Vert\theta
)\geq\frac{1}{\alpha}\log\int_{S}\frac{\theta^{\prime}}{\lambda^{\prime}%
}d\tilde{\nu}+\frac{1}{\alpha(\alpha-1)}\log\int_{\{\nu^{\prime}%
>0\}}\Big(\frac{\lambda^{\prime}}{\theta^{\prime}}\Big)^{\alpha-1}d\tilde{\nu
}.\label{23}%
\end{equation}

On the set $\{\nu^{\prime}>0\}=\{\lambda^{\prime}\theta^{\prime}>0\}$, define
\[
\varphi=\Big(\frac{\lambda^{\prime}}{\theta^{\prime}}\Big)^{\frac{\alpha
-1}{\alpha}},\qquad\psi=\Big(\frac{\theta^{\prime}}{\lambda^{\prime}%
}\Big)^{\frac{\alpha-1}{\alpha}},
\]
so that $\varphi\psi=1$ on $\{\nu^{\prime}>0\}$. By H\"{o}lder's inequality
with $1/p=1/\alpha$ and $1/q=(\alpha-1)/\alpha$, and with $p$ attached to
$\varphi$ and $q$ attached to $\psi$, we have
\begin{equation}
\int_{S}d\tilde{\nu}\leq\Big(\int_{\{\nu^{\prime}>0\}}\Big(\frac
{\lambda^{\prime}}{\theta^{\prime}}\Big)^{\alpha-1}d\tilde{\nu}\Big)^{\frac
{1}{\alpha}}\Big(\int_{S}\frac{\theta^{\prime}}{\lambda^{\prime}}d\tilde{\nu
}\Big)^{\frac{\alpha-1}{\alpha}}.\label{25}%
\end{equation}
Taking logs, dividing by $\alpha-1$ and using \eqref{23} gives that for any
$\theta\in\mathcal{P}$ with $\theta\gg\nu$,%
\[
\frac{1}{\alpha}\log\int_{S}e^{\alpha g}d\theta+R_{\alpha}(\nu\Vert\theta
)\geq\frac{1}{\alpha-1}\log\int_{S}d\tilde{\nu}=\frac{1}{\alpha-1}\log\int
_{S}e^{(\alpha-1)g}d\nu.
\]
If $\nu\not \ll \theta$ then $R_{\alpha}(\nu\Vert\theta)=\infty$, and again
the inequality holds.

Taking the infimum over all $\theta\in\mathcal{P}$ shows that the right hand
side of \eqref{11} is bounded below by the left hand side. Note that since $g$
is bounded $\lambda\{\nu^{\prime}>0\}=\lambda\{S\}$. Thus the choice
$\theta=\lambda$ gives equality in both \eqref{24} and \eqref{25}, hence in
\eqref{11}, and therefore identifies a minimizer.

Finally we show that the minimizer is unique. Assume that $\theta\gg\nu$
attains the infimum over $\mathcal{P}$. Then both \eqref{24} and \eqref{25}
must hold with equality. For \eqref{24} to hold with equality, $\theta\sim\nu$
must be true. Recall that H\"{o}lder's inequality will give an equality if and
only if $\theta^{\prime}/\lambda^{\prime}$ is constant on $\{\nu^{\prime
}>0\}=\{\lambda^{\prime}>0\}$. The only probability measure that satisfies
these conditions is $\theta=\lambda$, which shows that $\lambda$ attains the
infimum uniquely.

Next we consider (\ref{11}) for the same $\nu$, but for $\alpha\in(0,1)$. In
this case, we can no longer assume $\theta\gg\nu$. To show that the left hand
side of \eqref{11} is a lower bound for the right hand side, consider any
$\theta\in\mathcal{P}$. As with the case $\alpha>1$, let $\mu$ be a measure
dominating both $\nu$ and $\theta$, and define $\nu^{\prime},\theta^{\prime}$
and $\lambda^{\prime}$ with respect to this measure, where $d\lambda
=e^{-g}d\nu/Z$. Starting with the right hand side\ of \eqref{11},
\begin{align}
\log\int_{S}e^{\alpha g}d\theta &  \geq\log\int_{\{\nu^{\prime}\theta^{\prime
}>0\}}e^{\alpha g}d\theta\label{26}\\
&  =\log\int_{\{\nu^{\prime}\theta^{\prime}>0\}}\frac{1}{Z}\frac
{\theta^{\prime}}{\lambda^{\prime}}e^{(\alpha-1)g}d\nu,\nonumber
\end{align}
and
\begin{align*}
R_{\alpha}(\nu\Vert\theta)  &  =\frac{1}{\alpha(\alpha-1)}\log\int
_{\{\nu^{\prime}\theta^{\prime}>0\}}\Big(\frac{\nu^{\prime}}{\theta^{\prime}%
}\Big)^{\alpha}d\theta\\
&  =\frac{1}{\alpha(\alpha-1)}\log\int_{\{\nu^{\prime}\theta^{\prime}%
>0\}}Z^{\alpha-1}\Big(\frac{\lambda^{\prime}}{\theta^{\prime}}\Big)^{\alpha
-1}e^{(\alpha-1)g}d\nu.
\end{align*}
With $\tilde{\nu}$ again defined by $d\tilde{\nu}=e^{(\alpha-1)g}d\nu$,%
\begin{equation}
\frac{1}{\alpha}\log\int_{S}e^{\alpha g}d\theta+R_{\alpha}(\nu\Vert\theta
)\geq\frac{1}{\alpha}\log\int_{\{\nu^{\prime}\theta^{\prime}>0\}}\frac
{\theta^{\prime}}{\lambda^{\prime}}d\tilde{\nu}+\frac{1}{\alpha(\alpha-1)}%
\log\int_{\{\nu^{\prime}\theta^{\prime}>0\}}\Big(\frac{\lambda^{\prime}%
}{\theta^{\prime}}\Big)^{\alpha-1}d\tilde{\nu}. \label{27}%
\end{equation}
Define $\varphi=1$ and $\psi=(\lambda^{\prime}/\theta^{\prime})^{\alpha-1}$ on
the set $\{\nu^{\prime}\theta^{\prime}>0\}$. Using H\"{o}lder's inequality
with $p=1/\alpha$ attached to $\varphi$ and $q=1/(1-\alpha)$ attached to
$\psi$ gives
\[
\int_{\{\nu^{\prime}\theta^{\prime}>0\}}\Big(\frac{\lambda^{\prime}}%
{\theta^{\prime}}\Big)^{\alpha-1}d\tilde{\nu}\leq\Big(\int_{\{\nu^{\prime
}\theta^{\prime}>0\}}d\tilde{\nu}\Big)^{\alpha}\Big(\int_{\{\nu^{\prime}%
\theta^{\prime}>0\}}\frac{\theta^{\prime}}{\lambda^{\prime}}d\tilde{\nu
}\Big)^{1-\alpha}.
\]
Taking logs and dividing by $\alpha(\alpha-1)<0$ gives
\[
\frac{1}{\alpha(\alpha-1)}\log\int_{\{\nu^{\prime}\theta^{\prime}%
>0\}}\Big(\frac{\lambda^{\prime}}{\theta^{\prime}}\Big)^{\alpha-1}d\tilde{\nu
}\geq\frac{1}{\alpha-1}\log\int_{\{\nu^{\prime}\theta^{\prime}>0\}}d\tilde
{\nu}-\frac{1}{\alpha}\log\int_{\{\nu^{\prime}\theta^{\prime}>0\}}\frac
{\theta^{\prime}}{\lambda^{\prime}}d\tilde{\nu}.
\]
Using \eqref{27} gives
\begin{align}
\frac{1}{\alpha}\log\int_{S}e^{\alpha g}d\theta+R_{\alpha}(\nu\Vert\theta)  &
\geq\frac{1}{\alpha-1}\log\int_{\{\nu^{\prime}\theta^{\prime}>0\}}d\tilde{\nu
}\nonumber \\
&  =\frac{1}{\alpha-1}\log\int_{\{\nu^{\prime}\theta^{\prime}>0\}}%
e^{(\alpha-1)g}d\nu\nonumber \\
&  \geq\frac{1}{\alpha-1}\log\int_{S}e^{(\alpha-1)g}d\nu,\label{26a}
\end{align}
showing that \eqref{11} holds as an inequality. To show equality, substitute
$\lambda$ for $\theta$ and note that all the inequalities hold as equalities.

To show that $\lambda$ is the unique minimizer, note that any $\theta
\in\mathcal{P}$ satisfying all inequalities as equalities, must, in
particular, give equality in \eqref{26}, for which it is necessary that
$\theta \ll \nu$. Equality in \eqref{26a} implies $\nu \ll \theta$. 
For H\"{o}lder's inequality to hold with equality $\psi$
must be a constant, and the only probability measure satisfying these
conditions is $\lambda$. This completes the proof of (\ref{11}).

Toward proving (\ref{29}), note that (\ref{11}) implies
\[
\frac{1}{\alpha-1}\log\int_{S}e^{(\alpha-1)g}d\nu\leq\frac{1}{\alpha}\log
\int_{S}e^{\alpha g}d\theta+R_{\alpha}(\nu\Vert\theta),\qquad\nu,\theta
\in\mathcal{P},
\]
which is equivalent to
\[
\frac{1}{\alpha}\log\int_{S}e^{\alpha g}d\nu\geq\frac{1}{\alpha-1}\log\int
_{S}e^{(\alpha-1)g}d\theta-R_{\alpha}(\theta\Vert\nu),\qquad\nu,\theta
\in\mathcal{P}.
\]
Thus to prove part (\ref{29}), it suffices to show that the measure
$d\theta=e^{g}d\nu/Z$, and only this measure, gives equality in the above
display. The proof is similar to that of (\ref{11}), and therefore the details
are omitted. 

\bigskip

\noindent\textbf{Proof of Corollary \ref{cor:1}.} We give a proof of only the
rightmost inequality; the other inequality is proved analogously. First, if
$g$ is bounded the result follows from Theorem \ref{th1}. Otherwise, since the
claim holds trivially if the right hand side is infinite, assume it is finite.
Let $g^{M,N}=(g\vee-M)\wedge N$, for $M,N\geq0$. Then
\begin{equation}
\frac{1}{\alpha-1}\log\int_{S}e^{(\alpha-1)g^{M,N}}d\theta\leq\frac{1}{\alpha
}\log\int_{S}e^{\alpha g^{M,N}}d\nu+R_{\alpha}(\theta\Vert\nu).\label{05}%
\end{equation}
We first take $M\rightarrow\infty$ and use dominated convergence on both sides
of the inequality. To this end note that $g^{M,N}\leq g^{+}$. Since
$R_{\alpha}(\theta\Vert\nu)\geq0$ it must be true that $e^{\alpha g}$ is
$\nu$-integrable, and therefore so is $e^{\alpha g^{+}}$. Moreover, for
fixed $N$, $g^{M,N}\leq g^{0,N}$ and, using \eqref{05} with $M=0$, shows that
$e^{(\alpha-1)g^{0,N}}$ is $\theta$-integrable. As a result, \eqref{05} holds
with $g^{\infty,N}$ on both sides. Now we take $N\rightarrow\infty$ and use
monotone convergence (recall that $\alpha>1$). This gives the rightmost
inequality in \eqref{13} and completes the proof.

\bigskip

\noindent\textbf{Acknowledgment.} We are grateful to Ramon van Handel for
bringing to our attention the dualities in \cite{cha2} and \cite{mas2}.

\bigskip
\bibliographystyle{plain}
\bibliography{main}

\begin{thebibliography}{10}

\bibitem{ana2}
V.~Anantharam.
\newblock How large delays build up in a {$GI/G/1$} queue.
\newblock {\em Queueing Systems Theory Appl.}, 5(4):345--367, 1989.

\bibitem{asmgly}
S.~Asmussen and P.W. Glynn.
\newblock {\em Stochastic Simulation: Algorithms and Analysis}.
\newblock Applications of Mathematics. Springer Science+Business Media, LLC,
  2007.

\bibitem{bha}
A.~Bhattacharyya.
\newblock On some analogues of the amount of information and their use in
  statistical estimation.
\newblock {\em Sankhy\=a}, 8:1--14, 1946.

\bibitem{blagly}
J.H. Blanchet and P.~Glynn.
\newblock Efficient rare-event simulation for the maximum of heavy-tailed
  random walks.
\newblock {\em Ann. Appl. Prob.}, 18:1351--1378, 2008.

\bibitem{borsal}
A.~N. Borodin and P.~Salminen.
\newblock {\em Handbook of {B}rownian motion---facts and formulae}.
\newblock Probability and its Applications. Birkh\"auser Verlag, Basel, second
  edition, 2002.

\bibitem{cha2}
D.~Chafa{\"{\i}}.
\newblock Entropies, convexity, and functional inequalities: on
  {$\Phi$}-entropies and {$\Phi$}-{S}obolev inequalities.
\newblock {\em J. Math. Kyoto Univ.}, 44(2):325--363, 2004.

\bibitem{chodup}
K.~Chowdhary and P.~Dupuis.
\newblock Distinguishing and integrating aleatoric and epistemic variation in
  uncertainty quantification.
\newblock {\em ESAIM: Mathematical Modelling and Numerical Analysis},
  47:635--662, 2013.

\bibitem{demzei}
A.~Dembo and O.~Zeitouni.
\newblock {\em Large Deviations Techniques and Applications}.
\newblock Jones and Bartlett, Boston, 1993.

\bibitem{dupell4}
P.~Dupuis and R.~S. Ellis.
\newblock {\em A Weak Convergence Approach to the Theory of Large Deviations}.
\newblock John Wiley \& Sons, New York, 1997.

\bibitem{dupell}
P.~Dupuis and R.S. Ellis.
\newblock Large deviations for {M}arkov processes with discontinuous
  statistics, {II}: Random walks.
\newblock {\em Probab.\ Th.\ Rel.\ Fields}, 91:153--194, 1992.

\bibitem{dupjampet}
P.~Dupuis, M.~R. James, and I.~R. Petersen.
\newblock Robust properties of risk--sensitive control.
\newblock {\em Math.\ Control Signals Systems}, 13:318--332, 2000.

\bibitem{dupwan5}
P.~Dupuis and H.~Wang.
\newblock Subsolutions of an {I}saacs equation and efficient schemes for
  importance sampling.
\newblock {\em Math. Oper. Res.}, 32:1--35, 2007.

\bibitem{dvitod}
K.~Dvijotham and E.~Todorov.
\newblock A unified theory of linearly solvable optimal control.
\newblock {\em Artificial Intelligence (UAI)}, page~1, 2011.

\bibitem{frewen}
M.~I. Freidlin and A.~D. Wentzell.
\newblock {\em Random Perturbations of Dynamical Systems}.
\newblock Springer-Verlag, New York, 1984.

\bibitem{golpasyar}
L.~Golshani, E.~Pasha, and G.~Yari.
\newblock Some properties of {R}\'enyi entropy and {R}\'enyi entropy rate.
\newblock {\em Inform. Sci.}, 179(14):2426--2433, 2009.

\bibitem{kus15}
H.~J. Kushner.
\newblock Robustness and approximation of escape times and large deviations
  estimates for systems with small noise effects.
\newblock {\em SIAM J.\ Appl.\ Math.}, 44:160--182, 1984.

\bibitem{lievaj}
F.~Liese and I.~Vajda.
\newblock {\em Convex Statistical Distances}.
\newblock Teubner-Texte zur Mathematik. Teubner, 1987.

\bibitem{mas2}
P.~Massart.
\newblock {\em Concentration inequalities and model selection}, volume 1896 of
  {\em Lecture Notes in Mathematics}.
\newblock Springer, Berlin, 2007.
\newblock Lectures from the 33rd Summer School on Probability Theory held in
  Saint-Flour, July 6--23, 2003, With a foreword by Jean Picard.

\bibitem{nguwaijor}
X.~Nguyen, H.J. Wainwright, and M.I. Jordan.
\newblock Estimating divergence functionals and the likelihood ratio by convex
  risk minimization.
\newblock {\em IEEE Trans. Inform. Theory}, 56(11):5847--5861, 2010.

\bibitem{puhwhi}
A.A. Puhalskii and W.~Whitt.
\newblock Functional large deviation principles for first-passage-time
  processes.
\newblock {\em Ann. Appl. Probab.}, 7(2):362--381, 1997.

\bibitem{ren}
A.~R{\'e}nyi.
\newblock On measures of entropy and information.
\newblock In {\em Proc. 4th {B}erkeley {S}ympos. {M}ath. {S}tatist. and
  {P}rob., {V}ol. {I}}, pages 547--561, Berkeley, Calif., 1961. Univ.
  California Press.

\bibitem{rud2}
A.~Ruderman, M.~Reid, D.~Garc{\'\i}a-Garc{\'\i}a, and J.~Petterson.
\newblock Tighter variational representations of f-divergences via restriction
  to probability measures.
\newblock {\em arXiv preprint arXiv:1206.4664}, 2012.

\bibitem{shwwei}
A.~Shwartz and A.~Weiss.
\newblock {\em Large Deviations for Performance Analysis: Queues, Communication
  and Computing}.
\newblock Chapman and Hall, New York, 1995.

\bibitem{vaj}
I.~Vajda.
\newblock Distances and discrimination rates for stochastic processes.
\newblock {\em Stochastic Process. Appl.}, 35:47--57, 1990.

\bibitem{ervhar}
T.~van Erven and P.~Harremo{\"e}s.
\newblock R{\'e}nyi divergence and majorization.
\newblock In {\em Information Theory Proceedings (ISIT), 2010 IEEE
  International Symposium on}, pages 1335--1339. IEEE, 2010.

\bibitem{whi}
P.~Whittle.
\newblock Risk--sensitive linear/quadratic/guassian control.
\newblock {\em Adv. Appl. Prob.}, 13:777--784, 1981.

\end{thebibliography}

\end{document}